\documentstyle{article}
\begin{document}
\setlength{\baselineskip}{15pt}
\title{Simple evaluation of Casimir invariants in finite-dimensional Poisson 
   systems}
\author{Benito Hern\'{a}ndez--Bermejo \and V\'{\i}ctor Fair\'{e}n $^1$}
\date{}

\maketitle

{\em Departamento de F\'{\i}sica Fundamental, Universidad Nacional de 
Educaci\'{o}n a Distancia. Senda del Rey S/N, 28040 Madrid, Spain.}

\mbox{}

\begin{center} 
{\bf Abstract}
\end{center}
In this letter we present a procedure for the calculation of the Casimir 
functions of finite-dimensional Poisson systems which avoids the burden of 
solving a set of partial differential equations, as it is usually suggested 
in the literature. We show how a simple algebraic manipulation of the 
structure matrix reduces substantially the difficulty of the problem. 

\mbox{}

\mbox{}

\mbox{}


\mbox{}

\mbox{}

\mbox{}

\mbox{}

\mbox{}

\mbox{}

\mbox{}

\mbox{}

\mbox{}

\mbox{}

$^1$ To whom all correspondence should be addressed. E-mail vfairen@uned.es

\pagebreak
\begin{flushleft}
{\bf 1. Introduction}
\end{flushleft}

Poisson structures \cite{wei1,olv1} are ubiquitous in all fields of 
Mathematical Physics, from dynamical systems theory 
\cite{gyn1}--\cite{pla2} to fluid dynamics \cite{hs1,hs2}, 
magnetohydrodynamics \cite{hs2}--\cite{hyk3}, continuous media 
\cite{hyk3,dv1}, superconductivity \cite{hyk2}, superfluidity \cite{hyk1}, 
chromohydrodynamics \cite{ghk1}, etc. The association of a Poisson structure 
to a given physical problem (which is still an open question 
\cite{gyn1,hyg1,hoj2,pla1,byv1}) opens the possibility of obtaining a wide 
range of information about the system, which may be in the form of 
perturbative solutions \cite{lit1}--\cite{cyl1}, nonlinear stability analysis 
through the energy-Casimir method \cite{gyh1,ahyc1,hyct,hyw1}, 
bifurcation properties and characterization of chaotic behaviour 
\cite{dht1}--\cite{dht3}, or integrability results \cite{mag1,olv2}, to cite 
a few.

Mathematically, a finite-dimensional dynamical system is said to have a 
Poisson structure if it can be written in terms of a set of ODEs of the form: 
\begin{equation}
    \label{nham}
    \dot{x}^i = \sum_{j=1}^n J^{ij} \partial _j H \; , \;\:\; i = 1, \ldots , n, 
\end{equation} 
where $H(\mbox{{\bf x}})$, which is usually taken to be a time-independent 
first integral, plays the role of Hamiltonian function, and 
$J^{ij}(\mbox{{\bf x}})$ are the entries of a $n \times n$ skew-symmetric 
structure matrix ${\cal J}$ verifying the Jacobi 
equations:
\begin{equation}
     \label{jac}
     \sum_{l=1}^n ( J^{li} \partial_l J^{jk} + J^{lj} \partial_l J^{ki} + 
     J^{lk} \partial_l J^{ij} ) = 0 
\end{equation}
Here $ \partial_l $ means $ \partial / \partial x^l$ and indices $i,j,k$ run 
from 1 to $n$. Notice that, in particular, the rank of matrix ${\cal J}$ may 
not be maximum. For example, this is the case if the dimension of the 
system is odd, since the rank of a skew-symmetric matrix is always even. 
We shall denote in what follows the rank of matrix ${\cal J}$ by 
$2m$. It can be demonstrated \cite{olv1} that whenever the Poisson structure 
is singular (i.e., when $2m < n$) there exist $n-2m$ independent constants 
of motion known as Casimir (or distinguished) functions, which are present 
irrespective of the form of the Hamiltonian ---in other words, they are 
completely determined by the structure matrix. From an operational point of 
view, the Casimir functions $C(\mbox{{\bf x}})$ are the solutions of the 
set of partial differential equations ${\cal J} \cdot \nabla C = \mbox{{\bf 
0}}$, or equivalently \cite{olv1,hyg1,hoj1,litc}: 
\begin{equation}
    \label{caspde}
    \sum_{j=1}^n J^{ij} \partial _j C = 0 \; , \;\:\; i = 1, \ldots , 2m
\end{equation} 
Here we have assumed without loss of generality that the first $2m$ rows of 
${\cal J}$ are the linearly independent ones, a convention that we shall 
follow throughout. 

The characterization of the Casimir functions is of central importance in the 
analysis of Poisson structures. They do not only provide information about 
the structure of the solutions of the system (since they are first integrals, 
whose common level sets determine the symplectic foliation of the phase 
space). They constitute also the basis for establishing criteria for the 
nonlinear stability via the aforementioned energy-Casimir method; they allow 
the application of reduction of order procedures \cite{olv1,dhr1}; and they 
can be used in the determination of time-independent solutions of nonlinear 
field equations \cite{hhm1}. 

Resorting to system (\ref{caspde}) to obtain the Casimir functions is, in 
general, a rather inconvenient practice. We propose here a much simpler 
approach, which is developed from elementary linear algebraic considerations, 
which leads directly to a set of $n-2m$ {\em ordinary\/} differential 
equations. The application of our method, which is completely systematic, 
will be seen  to be {\em always\/} more efficient than the resolution of 
(\ref{caspde}). Moreover, the convenience of the procedure, when compared 
with the traditional approach, is greater for increasing dimension of the 
structure matrix.

\mbox{}

\begin{flushleft}
{\bf 2. Description of the method}
\end{flushleft}

Let us consider (\ref{nham}), and a region of the $n$-dimensional 
phase-space in which the rank of ${\cal J}$ is constant and equal to $2m<n$. 
If the $2m$ first rows of ${\cal J}$ are the linearly independent ones, then 
there exists a set of $2m \times (n-2m)$ functions 
$\gamma ^{i}_{k} (\mbox{{\bf x}})$, where $i = 2m+1 , \ldots , n$ and 
$k = 1 , \ldots , 2m$, such that
\begin{equation}
    \label{clj}
    J^{ij} = \sum_{k=1}^{2m} \gamma ^{i}_{k} J^{kj} \; , \;\:\; j = 1, \ldots , n
\end{equation} 
The importance of the proportionality functions $\gamma ^{i}_{k}$ was 
already noticed by Littlejohn \cite{litc}. Let us assume for the moment 
that they are known (their calculation is just a technical step for which we 
shall give a procedure later in this section). Then, the substitution of 
(\ref{clj}) into (\ref{nham}) gives immediately the following relations:
\begin{equation}
   \label{ecf1}
   \dot{x}^i = \sum_{k=1}^{2m} \gamma ^{i}_{k} \dot{x}^k 
              \; , \;\:\;\: i = 2m+1 , \ldots , n 
\end{equation}
These equations reveal the structure which is present in the system due 
to the fact that the rank of matrix ${\cal J}$ is not maximum, i.e., they 
express all interdependences among the system variables induced by the 
existence of the Casimir functions. We have therefore obtained a set of 
$(n-2m)$ ordinary differential equations for the Casimirs:
\begin{equation}
\label{ecf2}
   \mbox{d}{x}^i = \sum_{k=1}^{2m} \gamma ^{i}_{k} \mbox{d}{x}^k \; , 
   \;\:\;\: i = 2m+1 , \ldots , n 
\end{equation}
Note that each of these equations is to be integrated separately. 
It is not difficult to prove that (\ref{ecf2}) do lead to the Casimir 
functions: let $C^{(i)}(\mbox{\bf x})$ be a solution of the $i$-th equation, 
where $2m+1 \leq i \leq n$. Then there exists a function 
$\eta (\mbox{\bf x})$ such that:
\begin{equation}
   \label{grc}
    \mbox{d}C^{(i)} = \eta (\mbox{\bf x}) \left\{ \sum_{k=1}^{2m} 
    \gamma ^{i}_{k} \mbox{d}{x}^k - \mbox{d}{x}^i \right\}
\end{equation}
The $j$-th component of the vector ${\cal J} \cdot \nabla C^{(i)}$ will be:
\[
   ({\cal J} \cdot \nabla C^{(i)})^j = \sum_{k=1}^n J^{jk} \partial _k 
   C^{(i)} = \eta (\mbox{\bf x}) \left\{ \sum_{k=1}^{2m} J^{jk} 
   \gamma ^i_k - J^{ji} \right\} = 
\]
\begin{equation}
   \eta (\mbox{\bf x}) \left\{ J^{ij} - \sum_{k=1}^{2m} \gamma ^i_k J^{kj} 
   \right\} = 0 \;\: , \:\;\:\;\:\; \forall \; j = 1, \ldots ,n
\end{equation}
Here we have applied the original degeneracy relations (\ref{clj}). This 
demonstrates that the result of integrating each of the $n-2m$ equations 
(\ref{ecf2}) is one family of Casimir functions of matrix ${\cal J}$. We 
know, on the other hand, that there are $n-2m$ functionally independent 
Casimirs. From (\ref{grc}) it can be easily shown that the solutions of two 
different equations of the set (\ref{ecf2}) are always functionally 
independent. Consequently, the integration of equations (\ref{ecf2}) produces 
all the Casimirs of the system.

We end this section by indicating how functions $\gamma ^{i}_{k}$ can be   
calculated. To do so we proceed to write (\ref{clj}) in matrix form as: 
\begin{equation}
   \label{mclj}
   (\tilde{{\cal J}}_{2m})^T \cdot \Gamma = (\tilde{{\cal J}}_{n-2m})^T
\end{equation}
where $ \tilde{{\cal J}}_{2m}$ is the $2m \times n$ matrix composed by 
the first $2m$ rows of ${\cal J}$, $ \tilde{{\cal J}}_{n-2m}$ is the 
$(n-2m) \times n$ matrix composed by the last $(n-2m)$ rows of ${\cal J}$, 
and 
\begin{equation}
         \Gamma  = \left( \begin{array}{ccc}
                     \gamma ^{2m+1}_{1}  & \ldots  & \gamma ^{n}_{1} \\
                           \vdots        & \mbox{} & \vdots \\
                     \gamma ^{2m+1}_{2m} & \ldots  & \gamma ^{n}_{2m}
                    \end{array} \right)
\end{equation}
A rank analysis of the matrix equation (\ref{mclj}) shows immediately that 
there always exists a unique matrix $ \Gamma $ which is the solution. In 
fact, since $ \tilde{{\cal J}}_{2m}$ is a $2m \times n$ matrix, there are $(n-2m)^2$ 
redundant equations in (\ref{mclj}). If we assume again that these redundant  
equations are those corresponding to the last $(n-2m)$ rows of 
$( \tilde{{\cal J}}_{2m})^T$, we can write (\ref{mclj}) in the nonredundant 
form:
\begin{equation}
   \label{mclj2}
     ({\cal J}_{2m})^T \cdot \Gamma = ({\cal J}_{n-2m})^T
\end{equation}
where
\begin{equation}
   \label{rm}
     {\cal J}_{2m} = \left( \begin{array}{ccc}
                 J^{1,1}   & \ldots  & J^{1,2m}   \\
                   \vdots    & \mbox{} & \vdots       \\
                 J^{2m,1}  & \ldots  & J^{2m,2m}  
              \end{array} \right)  \; , \;\:\; 
     {\cal J}_{n-2m} =  \left( \begin{array}{ccc}
                 J^{2m+1,1}  & \ldots   & J^{2m+1,2m} \\
                    \vdots   & \mbox{}  & \vdots      \\
                 J^{n,1}     & \ldots   & J^{n,2m}    
              \end{array} \right)
\end{equation}
Since now ${\cal J}_{2m}$ is an invertible matrix, the solution is:
\begin{equation}
   \label{gamma}
     \Gamma =  ({\cal J}_{n-2m} \cdot {\cal J}_{2m}^{-1})^T 
\end{equation}

To summarize, our approach to the determination of the Casimir functions 
proceeds in two steps: {\em i)\/} calculation of $\Gamma$ through 
(\ref{gamma}); and {\em ii)\/} the integration of (\ref{ecf2}) ---each 
equation leading to an independent family of Casimirs. We shall now 
illustrate the procedure by means of some examples.

\mbox{}

\begin{flushleft}
{\bf 3. Examples}
\end{flushleft}
{\em (I) 3D Lotka-Volterra systems}

Nutku has demonstrated \cite{nut1} that the 3D Lotka-Volterra equations
\begin{eqnarray}
    \dot{x}^1 & = & x^1( \lambda + cx^2 + x^3) \nonumber   \\
    \dot{x}^2 & = & x^2( \mu     + x^1  + ax^3)             \\
    \dot{x}^3 & = & x^3( \nu     + bx^1 + x^2) \nonumber
\end{eqnarray}
are biHamiltonian when $abc=-1$ and $\nu = \mu b - \lambda a b$. In this 
case, the flow can be written as a Poisson system in two different ways:
\begin{equation}
    \dot{\mbox{\bf x}} = {\cal J}_1 \cdot \nabla H_1 = 
    {\cal J}_2 \cdot \nabla H_2 \;\: , 
\end{equation}
where:
\begin{equation}
     {\cal J}_1 = \left( \begin{array}{ccc}
                    0       & cx^1x^2 & bcx^1x^3  \\
                  -cx^1x^2  &   0     & -x^2x^3   \\
                 -bcx^1x^3  & x^2x^3  &  0    
              \end{array} \right) 
\end{equation}
\begin{equation}
     {\cal J}_2 = \left( \begin{array}{ccc}
                    0       & cx^1x^2(ax^3+ \mu) & cx^1x^3(x^2+ \nu)  \\
       -cx^1x^2(ax^3+ \mu)  &          0         &      x^1x^2x^3     \\
       -cx^1x^3(x^2+ \nu)   &    -x^1x^2x^3      &          0    
              \end{array} \right)  
\end{equation}
\begin{eqnarray}
      H_1 & = & abx^1+x^2-ax^3+ \nu \ln x^2 - \mu \ln x^3 \\
      H_2 & = & ab \ln x^1 -b \ln x^2 + \ln x^3 
\end{eqnarray}
Since the rank of both ${\cal J}_1$ and ${\cal J}_2$ is 2 everywhere in the 
positive orthant, there is always one independent Casimir. We shall apply 
our method to both Poisson structures. 

For ${\cal J}_1$ we have, by simple inspection:
\begin{equation}
   \mbox{(row3)} = \frac{x^3}{cx^1}\mbox{(row1)} + 
   \frac{bx^3}{x^2}\mbox{(row2)}
\end{equation}
In other words, $\gamma ^3_1 = x^3/cx^1$ and $\gamma ^3_2 = bx^3/x^2$. The 
equation we must solve is then:
\begin{equation}
   \mbox{d}x^3 = \frac{x^3}{cx^1}\mbox{d}x^1 + \frac{bx^3}{x^2}\mbox{d}x^2
\end{equation}
The integration of this equation is immediate and gives   
$ab \ln x^1 - b \ln x^2 + \ln x^3 = \mbox{constant}$, which is  
Nutku's result. Since any function of a Casimir is also a Casimir, the 
general solution would be:
\begin{equation}
   C = \Psi \left[ ab \ln x^1 - b \ln x^2 + \ln x^3 \right]
\end{equation}
with $\Psi$ a smooth one-variable function.

Similarly, for ${\cal J}_2$ we see that:
\begin{equation}
   \mbox{(row3)} = - \frac{x^3}{c(ax^3+ \mu)}\mbox{(row1)} + 
   \frac{x^3(x^2+ \nu)}{x^2(ax^3+ \mu)}\mbox{(row2)}
\end{equation}
Consequently, $\gamma ^3_1 = - x^3/(c(ax^3+ \mu))$ and $\gamma ^3_2 = 
x^3(x^2+ \nu)/(x^2(ax^3+ \mu))$. This implies that:
\begin{equation}
   \mbox{d}x^3 = - \frac{x^3}{c(ax^3+ \mu)} \mbox{d}x^1 + 
   \frac{x^3(x^2+ \nu)}{x^2(ax^3+ \mu)} \mbox{d}x^2
\end{equation}
After integration we arrive easily at $abx^1+x^2-ax^3+ \nu \ln x^2 - 
\mu \ln x^3 =$ constant, which is the solution. In general:
\begin{equation}
\label{c2b}
   C = \Psi \left[ abx^1+x^2-ax^3+ \nu \ln x^2 - \mu \ln x^3 \right]
\end{equation}

It is interesting to compare this procedure with the usual method of 
characteristics. We shall do it for ${\cal J}_2$. Since rank(${\cal J}_2$) 
is two in the domain of interest, the third equation of the system 
${\cal J}_2 \cdot \nabla C$ $=$ $0$ is a linear combination of the first and 
second ones, and can therefore be suppressed. The system of PDEs we have to 
solve in order to determine $C$ is then:
\begin{equation}
\label{edpe11}
   cx^1x^2(ax^3+ \mu) \frac{\partial C}{\partial x^2} + 
   cx^1x^3(x^2+ \nu) \frac{\partial C}{\partial x^3} = 0 
\end{equation}
\begin{equation}
\label{edpe12}
       -cx^1x^2(ax^3+ \mu) \frac{\partial C}{\partial x^1} + 
       x^1x^2x^3 \frac{\partial C}{\partial x^3} = 0 
\end{equation}
The characteristic equations for (\ref{edpe11}) are:
\begin{equation}
   \frac{\mbox{d}x^2}{cx^1x^2(ax^3+ \mu)} = 
   \frac{\mbox{d}x^3}{cx^1x^3(x^2+ \nu)} \:\; , \:\;\:
   \mbox{d}x^1 = 0
\end{equation}
Since $C$ is a function of three variables, we have to make two integrations 
from the characteristic equations. It can be found easily that 
$x^1=k_1$ and $x^2-ax^3+ \nu \ln x^2 - \mu \ln x^3 = k_2$, where $k_1$ and 
$k_2$ are constants of integration. Then, the general solution of equation 
(\ref{edpe11}) is of the form:
\begin{equation}
\label{c1e1}
    C^{(1)} = \Psi ^{(1)} [x^1 , x^2 - ax^3 + \nu \ln x^2 - \mu \ln x^3]
\end{equation}
Similarly, for the second PDE (\ref{edpe12}), the system of characteristic 
equations is:
\begin{equation}
 - \frac{\mbox{d}x^1}{cx^1x^2(ax^3+ \mu)} = 
   \frac{\mbox{d}x^3}{x^1x^2x^3} \:\; , \:\;\:
   \mbox{d}x^2 = 0
\end{equation}
We can obtain without difficulty that $x^2=k_1$ and $abx^1-ax^3 - \mu 
\ln x^3 = k_2$, and then the general solution of (\ref{edpe12}) is:
\begin{equation}
\label{c1e2}
    C^{(2)} = \Psi ^{(2)} [x^2 , abx^1 - ax^3 - \mu \ln x^3]
\end{equation}
Now we must take into account that the Casimir of the system is a 
simultaneous solution of both (\ref{edpe11}) and (\ref{edpe12}). This means 
that it must be a function of the $x^i$ complying to both formats 
(\ref{c1e1}) and (\ref{c1e2}). After inspection, one arrives directly to the 
solution (\ref{c2b}). We shall comment in Section 4 on the differences 
between both methods. 

\begin{flushleft}
{\em (II) A high-dimensional system: The light top}
\end{flushleft}

We shall now analyze in detail a six-dimensional example due to Weinstein 
\cite{wei1}: The equations of motion of a rigid body anchored at one point,   
which moves in a gravitational field. The system variables are the entries 
of the angular momentum in body coordinates, {\bf M} $= (M_1,M_2,M_3)$, as 
well as those of the gravitational force, also in body coordinates, {\bf F} 
$= (F_1,F_2,F_3)$. From now on, we will take the six variables in the 
following order: $(M_1,M_2,M_3,F_1,F_2,F_3)$. Then, the structure matrix and 
the Hamiltonian are, respectively:
\begin{equation}
     \label{je3}
     {\cal J} = \left( \begin{array}{cccccc}
                    0  &  M_3 & -M_2 &   0  &  F_3 & -F_2 \\
                  -M_3 &   0  &  M_1 & -F_3 &   0  &  F_1 \\
                   M_2 & -M_1 &   0  &  F_2 & -F_1 &   0  \\
                    0  &  F_3 & -F_2 &   0  &   0  &   0  \\
                  -F_3 &   0  &  F_1 &   0  &   0  &   0  \\
                   F_2 & -F_1 &   0  &   0  &   0  &   0
              \end{array} \right)  \;\;\: , 
\end{equation}
and 
\begin{equation}
   H = \sum_{i=1}^3 \left( \frac{M_i^2}{2I_i} + x_iF_i \right) 
\end{equation}
In $H$, the $I_i$ are the principal moments of inertia, and the $x_i$ are 
the coordinates of the body's center of mass measured from the anchor point 
(see \cite{wei1} and references therein for further details). 

We shall first apply our procedure for the determination of the Casimir 
functions of this system. For the sake of comparison, we shall later solve 
the same problem through the traditional method of characteristics.

\mbox{}

\noindent {\bf i)} Solution of the problem by the present method: 

Clearly, rank(${\cal J}$) $=4$, the third and the sixth rows being linear 
combinations of the rest. Then there are two independent Casimirs. We can 
find the $\gamma ^i_k$ by means of (\ref{gamma}):
\begin{equation}
    \Gamma = \left( \begin{array}{cccc}
      \gamma _1^3 & \gamma _2^3 & \gamma _4^3 & \gamma _5^3  \\
      \gamma _1^6 & \gamma _2^6 & \gamma _4^6 & \gamma _5^6  \\
      \end{array} \right) ^T = ({\cal J}_2 \cdot {\cal J}_4^{-1})^T
\end{equation}
where
\begin{equation}
     {\cal J}_4 = \left( \begin{array}{cccc}
                      0 &  M_3 &   0  &  F_3  \\
                   -M_3 &   0  & -F_3 &   0   \\ 
                      0 &  F_3 &   0  &   0   \\
                   -F_3 &   0  &   0  &   0
                         \end{array} \right) \; , \;\:
     {\cal J}_2 = \left( \begin{array}{cccc}
                    M_2 & -M_1 &  F_2 & -F_1  \\ 
                    F_2 & -F_1 &   0  &  0  
              \end{array} \right) 
\end{equation}
The solution is:
\begin{equation}
    \Gamma = \left( \begin{array}{cc}
            -F_1/F_3                &      0      \\
            -F_2/F_3                &      0      \\
            (F_1M_3-M_1F_3)/F_3^2   &  -F_1/F_3   \\
            (F_2M_3-M_2F_3)/F_3^2   &  -F_2/F_3   
             \end{array} \right) 
\end{equation}
We then have to solve independently the following two differential equations: 
\begin{equation}
 \label{eqcfm}
   \mbox{d}M_3 = - \frac{F_1}{F_3}\mbox{d}M_1 - \frac{F_2}{F_3}\mbox{d}M_2 + 
   \left( \frac{F_1M_3}{F_3^2}-\frac{M_1}{F_3} \right) \mbox{d}F_1 + 
   \left( \frac{F_2M_3}{F_3^2}-\frac{M_2}{F_3} \right) \mbox{d}F_2 
\end{equation}
\begin{equation}
 \label{eqcf2}
   \mbox{d}F_3 = - \frac{F_1}{F_3}\mbox{d}F_1 - \frac{F_2}{F_3}\mbox{d}F_2
\end{equation}
The last one is straightforward and gives a first Casimir: $C_1$ $=$ 
$F_1^2+F_2^2+F_3^2$ $=$ $\| \mbox{\bf F} \| ^2$. Now, if we expand 
(\ref{eqcfm}) and regroup terms we have: 
\begin{equation}
 \label{eqcfmm}
   F_1\mbox{d}M_1 + F_2\mbox{d}M_2 + F_3\mbox{d}M_3 + M_1\mbox{d}F_1 + 
   M_2\mbox{d}F_2 = M_3 \left( \frac{F_1}{F_3}\mbox{d}F_1 + 
   \frac{F_2}{F_3}\mbox{d}F_2\right)
\end{equation}
Making use of equation (\ref{eqcf2}) in the right-hand side of (\ref{eqcfmm}) 
leads immediately to $\mbox{d}(M_1F_1 + M_2F_2 + M_3F_3) = 0$. Thus, the 
second independent Casimir is $C_2$ $=$ $M_1F_1 + M_2F_2 + M_3F_3$ 
$=$ $\mbox{\bf M} \cdot \mbox{\bf F}$. We can write, as usual, the most 
general form of a Casimir as:
\begin{equation}
   C = \Psi [F_1^2+F_2^2+F_3^2,M_1F_1 + M_2F_2 + M_3F_3] \;\: ,
\end{equation}
where $\Psi$ is a smooth two-variable function. 

\mbox{}

\noindent {\bf ii)} Solution of the problem by the method of characteristics: 

We can now compare the previous procedure with the direct solution of the 
system of PDEs ${\cal J} \cdot \nabla C$ $=$ {\bf 0}. For this, we should 
begin by recalling the same observation than before: Since rank(${\cal J}$) 
$=4$, two of the equations of the system will be redundant ---which can be 
taken as those corresponding to the third and sixth rows of ${\cal J}$. 
Therefore, the system we have to solve is: 
\begin{equation}
\label{edpe1}
  M_3\frac{\partial C}{\partial M_2} -
  M_2\frac{\partial C}{\partial M_3} +
  F_3\frac{\partial C}{\partial F_2} -
  F_2\frac{\partial C}{\partial F_3} = 0 
\end{equation}
\begin{equation}
\label{edpe2}
- M_3\frac{\partial C}{\partial M_1} +
  M_1\frac{\partial C}{\partial M_3} -
  F_3\frac{\partial C}{\partial F_1} +
  F_1\frac{\partial C}{\partial F_3} = 0 
\end{equation}
\begin{equation}
\label{edpe3}
  F_3\frac{\partial C}{\partial M_2} -
  F_2\frac{\partial C}{\partial M_3} = 0 
\end{equation}
\begin{equation}
\label{edpe4}
- F_3\frac{\partial C}{\partial M_1} +
  F_1\frac{\partial C}{\partial M_3} = 0 
\end{equation}

The characteristic equations of (\ref{edpe1}) are:
\begin{equation}
\label{ece1}
   \frac{\mbox{d} M_2}{M_3} = - \frac{\mbox{d} M_3}{M_2} = 
   \frac{\mbox{d} F_2}{F_3} = - \frac{\mbox{d} F_3}{F_2} \;\; , \; 
   \mbox{d} M_1 = \mbox{d} F_1 = 0
\end{equation}
Since the unknown $C$ is a function of six variables, we have to find five 
constants from the characteristic equations (\ref{ece1}) in order to 
construct the general solution of the PDE (\ref{edpe1}). We immediately find 
from (\ref{ece1}) four of them:
\begin{equation}
   M_1=k_1 \;\; , \; F_1 = k_2 \;\; , \; M_2^2+M_3^2=k_3 \;\; , \; 
   F_2^2+F_3^2=k_4 
\end{equation}
We can derive a fifth one as follows:
\begin{eqnarray*}
   0 & \equiv & M_3 \mbox{d}F_3 - M_3 \mbox{d}F_3 + F_3 \mbox{d}M_3 - 
                F_3 \mbox{d}M_3 = \\
     &        & M_3 \mbox{d}F_3 + F_3 \mbox{d}M_3 + M_2 \mbox{d}F_2 + 
                F_2 \mbox{d}M_2 = \\
     &        & \mbox{d}(M_2F_2 + M_3F_3)
\end{eqnarray*}
Here we have made use of the characteristic equations (\ref{ece1}). The 
fifth constant is thus $k_5 = M_2F_2 + M_3F_3$. The general solution of the 
PDE (\ref{edpe1}) is then:
\begin{equation}
   C^{(1)} = \Psi ^{(1)}[M_1,F_1,M_2^2+M_3^2,F_2^2+F_3^2,M_2F_2 + M_3F_3]
\end{equation}

The second PDE (\ref{edpe2}) can be obtained from the first one (\ref{edpe1}) 
if we exchange the subindexes 1 and 2. Then we can directly write: 
\begin{equation}
   C^{(2)} = \Psi ^{(2)}[M_2,F_2,M_1^2+M_3^2,F_1^2+F_3^2,M_1F_1 + M_3F_3]
\end{equation}

For the third equation (\ref{edpe3}) we now have:
\begin{equation}
\label{ece3}
   \frac{\mbox{d} M_2}{F_3} = - \frac{\mbox{d} M_3}{F_2} \;\; , \; 
   \mbox{d} M_1 = \mbox{d} F_1 = \mbox{d} F_2 = \mbox{d} F_3 = 0
\end{equation}
This leads to:
\begin{equation}
    M_1 = k_1 \;\; , \; F_1 = k_2 \;\; , \;  F_2 = k_3 \;\; , \; F_3 = k_4
\end{equation}
Since $F_2$ and $F_3$ are constants, we also arrive at $k_5 = M_2F_2 + 
M_3F_3$. Consequently, the general solution of the PDE (\ref{edpe3}) is:
\begin{equation}
   C^{(3)} = \Psi ^{(3)}[M_1,F_1,F_2,F_3,M_2F_2 + M_3F_3]
\end{equation}

And finally, we again obtain the fourth PDE (\ref{edpe4}) from the third one 
(\ref{edpe3}) by permutation of the subindexes 1 and 2. Therefore: 
\begin{equation}
   C^{(4)} = \Psi ^{(4)}[M_2,F_1,F_2,F_3,M_1F_1 + M_3F_3]
\end{equation}

Now, the Casimir functions are simultaneous solutions of {\em all\/} the PDEs 
(\ref{edpe1}-\ref{edpe4}). Then, we now have to compare the four solutions 
$C^{(i)}$, for $i = 1, \ldots ,4$, and look for those functions of {\bf M} 
and {\bf F} compatible with all of them. After inspection, it is not 
difficult to arrive to the two most ovbious possibilities: 
$\| \mbox{\bf F} \| ^2$ and $\mbox{\bf M} \cdot \mbox{\bf F}$, which are the 
two independent Casimirs already known.

\begin{flushleft}
{\bf 4. Final remarks}
\end{flushleft}

We have seen how our algebraic approach allows the calculation of the 
Casimir functions in a quite natural and rapid way. In fact, we believe that 
this procedure gives some insight on how a symplectic foliation arises from 
the degereracy present in a singular Poisson structure. 

A comparison with the traditional method relying on the system of PDEs 
(\ref{caspde}) seems to be convenient. If we wish to solve equations 
(\ref{caspde}), the two simplest strategies are separation of variables and 
the method of characteristics. 

Separation of variables, which is rather lengthy even for simple PDEs and 
usually requires an eigenvalue analysis of the resulting ODEs, is clearly 
much less efficient than our procedure. 

On the other hand, we have already given in the examples a comparative 
solution of the problems by both our approach and the method of 
characteristics. Before entering in more quantitative and general arguments, 
two observations can be drawn from the examples: The first one is that our 
method is clearly less computation consuming than that of the 
characteristics. Notice that our technique reduces the problem to the 
solution of one ODE per Casimir. The number of ODEs which has been necessary 
to handle and the number of quadratures which must be found by the method of 
characteristics is certainly higher, in both examples. The second important 
remark is that both techniques do not lead to the same set of equations, 
i.e., our method is {\em not \/} a shortcut for the obtainment of the 
characteristic equations, as it can be easily checked. 

Let us compare in a quantitative way the complexity of both methods. We 
shall give as a measure of such complexity the number of quadratures which 
have to be calculated in every case to determine the solution. This number 
is $N_{a} = n-2m$ for our algebraic method, namely the corank of the 
structure matrix, as we already know. 

In the method of the characteristics, on the other hand, we have to solve 
system (\ref{caspde}), which consists of $2m$ nonredundant PDEs (the 
remaining $n-2m$ equations are redundant due to the degereracy in rank of the 
structure matrix, and can therefore be suppressed, as we have seen in the 
examples). In order to compute the total number of quadratures in the method 
of characteristics, let us consider the $i$-th PDE of system (\ref{caspde}). 
Its characteristic equations are of the form:
\begin{equation}
   \frac{\mbox{d} x^1}{J^{i1}} = \ldots = \frac{\mbox{d} x^{i-1}}{J^{i,i-1}} 
   = \frac{\mbox{d} x^{i+1}}{J^{i,i+1}} = \ldots = 
   \frac{\mbox{d} x^n}{J^{in}} \; , \;\; \mbox{d} x^i=0
\end{equation} 
Since $C$ is a function of $n$ variables, we need $n-1$ quadratures. However, 
we always have a trivial one, which is $x^i =$ constant. Therefore, we only 
have to carry out $n-2$ quadratures per PDE, in general. Consequently, the 
total number of quadratures is $N_{c} = 2m(n-2)$ for the method of 
characteristics. It is then straightforward to verify that 
\begin{equation}
\label{comp}
    \frac{N_a}{N_c} < 1
\end{equation}
in all nontrivial cases (the only situation in which (\ref{comp}) is not 
satisfied for a singular Poisson structure, is the unimportant case 
corresponding to a null structure matrix). 
When the number of Casimirs is large, for example 
if $2m = 2$, we obtain $N_a/N_c = 1/2$. When such a number is medium, 
i.e. for $2m \simeq n/2$, we have that $N_a/N_c \simeq 1/(n-2)$, thus 
decreasing with increasing size of the structure matrix. Finally, when the 
number of Casimirs is small, say $2m \simeq (n-1)$, we arrive at $N_a/N_c 
\simeq 1/[(n-1)(n-2)]$. In this case the ratio decreases as $n^{-2}$ as 
$n$ grows, and our approach is now much more economic for a large 
structure matrix. 

\mbox{}

\mbox{}

\begin{flushleft}
{\bf Acknowledgements}
\end{flushleft}

This work has been supported by the DGICYT (Spain), under grant PB94-0390, 
and by the EU Esprit WG 24490 (CATHODE-2). B. H. acknowledges a doctoral 
fellowship from Comunidad Aut\'{o}noma de Madrid. The authors also wish to 
acknowledge an anonymous referee for useful suggestions and comments. 

\pagebreak


\begin{thebibliography}{99}
   \bibitem{wei1} A. Weinstein, J. Diff. Geom. 18 (1983) 523.
   \bibitem{olv1} P. J. Olver, Applications of Lie Groups to Differential 
      Equations, 2nd Ed. (Springer-Verlag, New York, 1993).
   \bibitem{gyn1} H. G\"{u}mral and Y. Nutku, J. Math. Phys. 34 (1993) 5691.
   \bibitem{hyg1} F. Haas and J. Goedert, Phys. Lett. A 199 (1995) 173.
   \bibitem{gyh1} J. Goedert, F. Haas, D. Hua, M. R. Feix and L. Cair\'{o}, 
      J. Phys. A: Math. Gen. 27 (1994) 6495.
   \bibitem{hoj1} S. A. Hojman, J. Phys. A: Math. Gen. 24 (1991) L249.
   \bibitem{hoj2} S. A. Hojman, J. Phys. A: Math. Gen. 29 (1996) 667.
   \bibitem{pla1} M. Plank, J. Math. Phys. 36 (1995) 3520.
   \bibitem{pla2} M. Plank, Nonlinearity 9 (1996) 887.
   \bibitem{hs1} D. D. Holm, Physica D 17 (1985) 1.
   \bibitem{hs2} D. D. Holm, Phys. Lett. A 114 (1986) 137.
   \bibitem{myg1} P. J. Morrison and J. M. Greene, Phys. Rev. Lett. 45 
      (1980) 790.
   \bibitem{hhm1} R. D. Hazeltine, D. D. Holm and P. J. Morrison, J. Plasma 
      Physics 34 (1985) 103.
   \bibitem{hyk3} D. D. Holm and B. A. Kupershmidt, Physica D 6 (1983) 347.
   \bibitem{dv1} I. E. Dzyaloshinskii and G. E. Volovick, Ann. Physics 125 
      (1980) 67.
   \bibitem{hyk2} D. D. Holm and B. A. Kupershmidt, Phys. Lett. A 93 (1983) 177.
   \bibitem{hyk1} D. D. Holm and B. A. Kupershmidt, Phys. Lett. A 91 (1982) 425.
   \bibitem{ghk1} J. Gibbons, D. D. Holm and B. Kupershmidt, Phys. Lett. A 
      90 (1982) 281.
   \bibitem{byv1} B. Hern\'{a}ndez--Bermejo and V. Fair\'{e}n, Phys. Lett. 
      A 234 (1997) 35.
   \bibitem{lit1} R. G. Littlejohn, J. Math. Phys. 20 (1979) 2445.
   \bibitem{lit2} R. G. Littlejohn, J. Math. Phys. 23 (1982) 742.
   \bibitem{cyl1} J. R. Cary and R. G. Littlejohn, Ann. Physics 151 (1983) 1.
   \bibitem{ahyc1} H. D. I. Abarbanel, D. D. Holm, J. E. Marsden and T. 
      Ratiu, Phys. Rev. Lett. 52 (1984) 2352.
   \bibitem{hyct} D. D. Holm, J. E. Marsden, T. Ratiu and A. Weinstein, Phys. 
      Rep. 123 (1985) 1.
   \bibitem{hyw1} D. D. Holm and K. B. Wolf, Physica D 51 (1991) 189.
   \bibitem{dht1} D. David, D. D. Holm and M. V. Tratnik, Phys. Lett. A 137 
      (1989) 355.
   \bibitem{dht2} D. David, D. D. Holm and M. V. Tratnik, Phys. Lett. A 138  
      (1989) 29.
   \bibitem{dht3} D. David, D. D. Holm and M. V. Tratnik, Phys. Rep. 187  
      (1990) 281.
   \bibitem{mag1} F. Magri, J. Math. Phys. 19 (1978) 1156.
   \bibitem{olv2} P. J. Olver, Phys. Lett. A 148 (1990) 177.
   \bibitem{litc} R. G. Littlejohn, AIP Conf. Proc. 88 (1982) 47.
   \bibitem{dhr1} D. David and D. D. Holm, J. Nonlinear Sci. 2 (1992) 241.
   \bibitem{nut1} Y. Nutku, Phys. Lett. A 145 (1990) 27. 
\end{thebibliography}
\end{document}